\documentclass[12pt, letterpaper]{article}
\usepackage{amsfonts}
\usepackage{amssymb}
\usepackage{latexsym}  

\newcommand\BB{{\mathcal B}}
\newcommand\PP{{\mathcal P}}
\newcommand\TT{{\mathcal T}}

\newcommand\liminv{\overleftarrow{\rule{0cm}{3.7mm}\lim}}  
\newcommand\onto{\twoheadrightarrow}  

\newcommand\cl{\mathrm{cl}} 
\newcommand\intr{\mathrm{int}} 
\newcommand\cf{\mathrm{cf}}   
\newcommand\st{\mathrm{st}}  

\newcommand\iv{^{-1}} 

\newcommand\eop{$\ \ {\vcenter
   {\hrule
   \hbox{\vrule height 9pt \kern 9pt \vrule height 9pt}
   \hrule}}$\vskip 1.0 pt}

{\end{enumerate}}
\newenvironment{itemizz}{\begin{itemize}\setlength{\itemsep}{-1mm}} %
{\end{itemize}}                              

\newtheorem{theorem}{Theorem} 
\newtheorem{definition}[theorem]{Definition}
\newtheorem{lemma}[theorem]{Lemma}

\newenvironment{proof}{{\bf Proof.}}{\eop\medskip}
\newenvironment{proofof}[1]{\medskip \textbf{Proof of #1.}}{\eop\medskip}

\begin{document}

\title{Small Locally Compact Linearly Lindel\"of Spaces\footnote{  %
2000 Mathematics Subject Classification:   %
Primary 54D20, 54D80; Secondary  03E55.    %
Key Words and Phrases: linearly Lindel\"of, weak P-point, J\'onsson cardinal.}
}

\author{
Kenneth Kunen\footnote{University of Wisconsin,  Madison, WI  53706, U.S.A.,
\ \ kunen@math.wisc.edu}
\thanks{Author partially supported by NSF Grant DMS-0097881.}
}

\maketitle

\begin{abstract}
There is a locally compact Hausdorff space
of weight $\aleph_\omega$
which is
linearly Lindel\"of and not Lindel\"of.
\end{abstract}

\noindent
We shall prove:

\begin{theorem}
\label{thm-main}
There is a compact Hausdorff space $X$ and a point $p$ in $X$ such that:
\begin{itemizz}
\item[1.] $\chi(p,X) = w(X) = \aleph_\omega$.
\item[2.] For all regular $\kappa > \omega$, no $\kappa$-sequence
of points distinct from $p$ converges to $p$.
\end{itemizz}
\end{theorem}

As usual, $\chi(p,X)$, the \textit{character} of $p$ in $X$,
is the least size of a local base at $p$,
and $w(X)$, the \textit{weight} of $X$,
is the least size of a base for $X$.
This theorem with ``$\beth_\omega$'' replacing ``$\aleph_\omega$''
was proved in \cite{KUN2}.
Arhangel'skii and  Buzyakova \cite{AB} point out that if $X,p$
satisfy (2) of the theorem, then the space $X \backslash \{p\}$ is
linearly Lindel\"of  and locally compact; if in addition 
$\chi(p,X) > \aleph_0$, then $X \backslash \{p\}$ is not 
Lindel\"of.  (2) requires $\cf(\chi(p,X)) = \omega$, because there must be
a sequence of type $\cf(\chi(p,X))$ converging to $p$.
Thus, in (1) of the theorem, 
$\aleph_\omega$ is the smallest possible uncountable value for
$\chi(p,X)$ and $w(X)$.

As in \cite{KUN2}, the $X$ of the theorem will
be constructed as an inverse limit, using the following terminology:

\begin{definition}
An \emph{inverse system} is a sequence 
$\langle X_n, \pi^{n+1}_n : n \in \omega \rangle$,
where each $X_n$ is a compact Hausdorff space, and
each $\pi^{n+1}_n$ is a continuous map from $X_{n+1}$ onto $X_n$.
\end{definition}

Such an inverse systems yields a compact Hausdorff space, 
$X _\omega = \liminv_n X_n $,
and maps $\pi^\omega_m : X_\omega \onto X_m$
for $m < \omega$ and $\pi^n_m : X_n \onto X_m$ for $m \le n < \omega$.
Exactly as in \cite{KUN2}, one easily proves:

\begin{lemma}
\label{lemma-properties}
Suppose that
$\langle X_n, \pi^{n+1}_n : n \in \omega \rangle$
is an inverse system and $p \in X = X_\omega$, with the
$p_n = \pi^\omega_n (p) \in X_n$ satisfying:
\begin{itemizz}
\item[A.] Each $p_n$ is a weak $P_{\aleph_{n}}$-point in $X_n$.
\item[B.] Each $w(X_n) < \aleph_\omega$.
\item[C.] Each $(\pi^n_0)\iv\{p_0\}$ is nowhere dense in $X_n$.
\end{itemizz}
Then $X,p$ satisfies Theorem \ref{thm-main}.
\end{lemma}

As usual, $y \in Y$ is a \textit{weak} $P_\kappa$-point iff $y$ is
not in the closure of any subset of $Y\backslash\{y\}$ of size
less than $\kappa$, and
$y$ is a  $P_\kappa$-point iff the intersection of fewer than $\kappa$
neighborhoods of $y$ is always a neighborhood of $y$.
These properties are trivial for $\kappa = \aleph_0$.  The
terms ``$P$-point'' and ``weak $P$-point'' denote
``$P_{\aleph_1}$-point'' and ``weak $P_{\aleph_1}$-point'', 
respectively.

Every $P_\kappa$-point  is a weak $P_\kappa$-point, but
as pointed out in \cite{KUN2}, one cannot have 
each $p_n$ being a $P_{\aleph_{n}}$-point, as that would contradict (C).
In the construction we describe, it will be natural to make every
$p_n$ fail to be a $P$-point in $X_n$.

We shall build the $X_n$ and $p_n$ inductively using the following:

\begin{lemma}
\label{lemma-P-cover}
Assume that $y \in F \subseteq Y$, where $Y$ is compact Hausdorff,
$w(Y) \le \aleph_n$, and $\intr(F) = \emptyset$.  Then 
there is a compact Hausdorff space $X$, a point $x \in X$,
and a continuous $g: X \to Y$ such that:
\begin{itemizz}
\item[1.] $g(X) = Y$ and $g(x) = y$.
\item[2.] $g\iv(F)$ is nowhere dense in $X$.
\item[3.] $w(X) = \aleph_n$.
\item[4.] In $X$, $x$ is a weak $P_{\aleph_n}$-point and not a $P$-point.
\end{itemizz}
\end{lemma}

\begin{proofof}{Theorem \ref{thm-main}}
Inductively construct
the inverse system as in Lemma \ref{lemma-properties},
with each  $w(X_n) = \aleph_n$.
$X_0$ can be the Cantor set. When $n > 0$ and we are
given $X_{n-1},p_{n-1}$, we apply Lemma  \ref{lemma-P-cover} with
$F = (\pi^{n-1}_0)\iv\{p_0\}$.
\end{proofof}

Of course, we still need to prove Lemma \ref{lemma-P-cover}.
We remark that we do not assume that $F$ is closed,
although that was true in our proof of Theorem \ref{thm-main}.
Even if $F$ is dense in $Y$ in
Lemma \ref{lemma-P-cover},
we still get (2) --- that is $\intr (\cl (g\iv(F))) = \emptyset$.

When $n = 0$ in Lemma \ref{lemma-P-cover},
the ``weak $P_{\aleph_0}$-point'' is trivial,
and the lemma is easily proved by an Aleksandrov duplicate construction.
A more convoluted proof is:
Let $D \subseteq Y \backslash F$ be dense in $Y$ and countable.
Let $g$ map $\omega$ onto $D$ and extend $g$ to a map
$\beta g : \beta\omega \onto Y$.  Choosing $x$ to be any
point in $(\beta g)\iv(\{y\})$ yields (1)(2)(4),
but $\beta\omega$ has weight  $2^{\aleph_0}$.
Now, we can take a countable elementary submodel of the
whole construction to get an $X$ of weight $\aleph_0$.
Our proof for a general $n$ will follow this pattern.

As usual, $\beta\kappa$ denotes the \v Cech compactification
of a discrete $\kappa$, and $\kappa^* = \beta\kappa \backslash \kappa$.
Equivalently, $\beta\kappa$ is
the space of ultrafilters on $\kappa$ and $\kappa^*$ is the space
of nonprincipal ultrafilters.
If $g : \kappa \to Y$, where $Y$ is compact Hausdorff, then 
$\beta g$ denotes the unique extension of $g$ to a continuous
map from $\beta \kappa$ to $Y$.
Our weak $P_\kappa$-points in Lemma \ref{lemma-P-cover} will
be \textit{good} ultrafilters in the sense of Keisler \cite{KEISLER1}:

\begin{definition}
An ultrafilter $x$ on $\kappa$ is \emph{good} iff
for all $H: [\kappa]^{< \omega} \to x$, there is a
$K : \kappa \to x$ such that for each 
$s = \{\alpha_1 , \ldots , \alpha_n\} \in [\kappa]^{< \omega}$,
$K(\alpha_1) \cap \cdots \cap K(\alpha_n) \subseteq H(s)$.
\end{definition}

The following is well-known.

\begin{lemma}
Let $\kappa$ be any infinite cardinal.
\begin{itemizz}
\item[1.] There are ultrafilters $x$ on $\kappa$ which are
both good and countably incomplete.
\item[2.] Any $x$ as in (1) is a weak $P_\kappa$ point and
not a $P$-point in $\beta\kappa$.
\end{itemizz}
\end{lemma}

In (2), $x$ is not a $P$-point by countable incompleteness,
and proofs that it is a weak $P_\kappa$ point can be found in
\cite{BK1, BK2, DOW}.  For (1), see \cite{CK}, Theorem 6.1.4;
also, \cite{BK1, BK2} construct good ultrafilters with various
additional properties.

We first point out (Lemma \ref{lemma-reg}) that
taking $x$ to be a good ultrafilter
on $\omega_n$ will give us (1)(2)(4) of Lemma \ref{lemma-P-cover}.
Unfortunately, $w(\beta\omega_n) = 2^{\aleph_n}$,
so we shall take an elementary submodel to 
bring the weight down.  Omitting the elementary submodel,
our argument is as in
\cite{KUN2}, which obtained the $X$ of Theorem \ref{thm-main}
with $w(X) = \beth_\omega$, rather than $\aleph_\omega$.
A related use of elementary submodels to reduce
the weight occurs in \cite{JK}.

Before we consider the weight problem, we explain how
to map the good ultrafilter onto the given point $y$.
This part of the argument works for any regular ultrafilter.

\begin{definition}
\label{def-reg}
An ultrafilter $x$ on $\kappa$ is \emph{regular}
iff there are $E_\alpha\in x$ for $\alpha < \kappa$ such that
$\{\alpha : \xi \in E_\alpha\}$ is finite for all $\xi < \kappa$.
\end{definition}

Such an $x$ is countably incomplete
because $\bigcap_{n < \omega} E_n = \emptyset$.
For the following, see
Exercise 6.1.3 of \cite{CK} or
the proof of Lemma 2.1 in Keisler \cite{KEISLER2}:

\begin{lemma}
\label{lemma-good-reg}
If $x$ is a countably incomplete good ultrafilter on $\kappa$,
then $x$ is regular.
\end{lemma}

\begin{lemma}
\label{lemma-reg}
Let $x$ be a regular ultrafilter on $\kappa$.  
Assume that $y \in F \subseteq Y$, where $Y$ is compact Hausdorff,
$w(Y) \le \kappa$, and $\intr(F) = \emptyset$.  Then 
there is a map $g : \kappa \to Y$ such that
\begin{itemizz}
\item[A.] $\beta g$ maps $\beta \kappa$ onto $Y$.
\item[B.] $(\beta g)(x) = y$.
\item[C.] $g(\xi) \notin F$ for all $\xi \in \kappa$.
\item[D.] $g\iv(F)$ is nowhere dense in $\beta\kappa$.
\end{itemizz}
\end{lemma}
\begin{proof}
Of course, (D) follows from (C) because $g\iv(F) \subseteq \kappa^*$.
Fix $A \subseteq \kappa$ with $A \notin x$ and $|A| = \kappa$.
Let $\{E_\alpha : \alpha < \kappa\}$ be as in Definition \ref{def-reg},
with each
$E_\alpha \cap A = \emptyset$.  Let $\{U_\alpha : \alpha < \kappa\}$
be an open base at $y$ in $Y$.
Let $D \subseteq Y \backslash F$ be dense in $Y$.
Choose $g : \kappa \to Y$ such that $g$ maps $A$ onto $D$
(ensuring (A)) and each
$g(\xi) \in \bigcap\{U_\alpha : \xi \in E_\alpha\} \setminus F$
(ensuring (B)(C)).
\end{proof}

To apply the elementary submodel technique (as in Dow \cite{DOW3}),
we put the construction of Lemma \ref{lemma-reg}
inside an $H(\theta)$, where $\theta$ is a suitably large
regular cardinal.  Let $M \prec H(\theta)$, with
$\kappa \subset M$ and $|M| = \kappa$, such that
$M$ contains $Y$ and its topology $\TT$, along with $F,g,x,y$.
Let $\BB = \PP(\kappa) \cap M$, let  $\st(\BB)$ denote its Stone space,
and let $\Gamma: \beta \kappa  \onto \st(\BB)$ be the natural map; so
$\Gamma(x) = x \cap \BB = x \cap M$.
Since $\TT \cap M$ is a base for $Y$ (by $w(Y) \le \kappa$), we have
$\Gamma(z_1) = \Gamma(z_2) \to (\beta g)(z_1) = (\beta g)(z_2)$,
so that $\beta g$ yields
a map $\widetilde g : \st(\BB) \to Y$ with 
$\beta g = \widetilde g  \circ \Gamma$.  Note that $\BB$ contains all
finite subsets of $\kappa$, so that
$\st(\BB)$ is some compactification of a discrete $\kappa$.
It is easily seen that we 
still have (A--D), replacing $\beta g$ by  $\widetilde g$,
$\beta \kappa$ by $\st(\BB)$, and $x$ by $\Gamma(x)$.
Note that $\Gamma(x)$ must be countably incomplete by
$M \prec H(\theta)$, so that $\Gamma(x)$ will not be a $P$-point
in $\st(\BB)$.
But to prove Lemma \ref{lemma-P-cover} (letting $\kappa = \aleph_n$),
we also need $\Gamma(x)$ to be a weak $P_\kappa$-point in $\st(\BB)$.
We may assume that $x \in \beta\kappa$ is good, 
so it is a weak $P_\kappa$-point there.
But we need to show that in $\st(\BB)$, $\Gamma(x)$
is not a limit point of any set of size $\lambda < \kappa$.
Our argument here needs to assume
that $M$  is $\lambda$-covering and that $\lambda^+$ is not
a J\'onsson cardinal.  These two assumptions will cause
no problem when $\lambda < \aleph_\omega$.

As usual, $M \prec H(\theta)$ is \textit{$\lambda$-covering}
iff for all $E \in [M]^\lambda$, there is an $F \in [M]^\lambda$
such that $E \subseteq F$ and $F \in M$.
By taking a union of an elementary chain of type $\lambda^+$
(see \cite{DOW3}, \S3), we see that there is an
$M \prec H(\theta)$ with $|M| = \lambda^+$ such that
$M$ is $\lambda$-covering.

$\kappa$ is called a \textit{J\'onsson cardinal} iff
for all $\psi : [\kappa]^{< \omega} \to \kappa$, there
is a $W \in [\kappa]^\kappa$ such that
$\psi([W]^{< \omega})$ is a proper subset of $\kappa$.
By Tryba \cite{TR} (or see \cite{KAN}):

\begin{lemma}
\label{lemma-jon}
No successor to a regular cardinal is J\'onsson.
\end{lemma}

In particular, each $\aleph_n$ is not a J\'onsson cardinal;
this fact is much older and is easily proved by induction on $n$.

\begin{lemma}
\label{lemma-reflect}
Let $\kappa$ be infinite and  $x \in \beta\kappa$ a good ultrafilter
on $\kappa$.  Fix an infinite $\lambda < \kappa$ and let
$\theta > 2^\kappa$ be regular.  Let $M \prec H(\theta)$, with
$x, \kappa \in M$ and $\kappa \subset M$.  Assume
that $M$ is $\lambda$-covering and $\lambda^+$ is not
a J\'onsson cardinal.
Let $\BB = \PP(\kappa) \cap M$, and let 
$\Gamma: \beta \kappa  \onto \st(\BB)$ be the natural map.
Then $\Gamma(x)$ is a weak $P_{\lambda^+}$-point of $\st(\BB)$.
\end{lemma}
\begin{proof}
Fix $S \subseteq \st(\BB)\backslash\{\Gamma(x)\}$ with 
$|S| \le \lambda$.  We shall show that $\Gamma(x)$ is not in 
the closure of $S$.
For each $z \in S$, choose $F_z \in \Gamma(x) = x \cap \BB = x \cap M$
such that $F_z \notin z$.
Since $M$ is  $\lambda$-covering, we can get
$\langle G_\xi : \xi < \lambda \rangle \in M$ such that
each $G_\xi \in x$ and
$\forall z \in S\, \exists \xi < \lambda\, [G_\xi = F_z]$.
Since  $\lambda^+$ is not J\'onsson and $\lambda^+ \in M$, we 
can fix $\psi \in M$ such that
$\psi : [\lambda^+]^{< \omega} \to \lambda$ and  such that
$\psi([W]^{< \omega}) = \lambda$ for all $W \in [\lambda^+]^{\lambda^+}$.
Define $H(s) = G_{\psi(s)}$.  Then $H \in M$ and
$H : [\lambda^+]^{< \omega} \to \Gamma(x)$.
Since $x$ is good, we can find
$\langle K_\alpha : \alpha < \lambda^+ \rangle \in M$ such 
that each $K_\alpha$ is in $x$ (and hence in $\Gamma(x) = x \cap M$),
and such that $K_{\alpha_1} \cap \cdots \cap K_{\alpha_n}  \subseteq
H(\{\alpha_1,  \ldots , \alpha_n \})$ for each $n$ and each
$\alpha_1 , \ldots , \alpha_n \in \lambda^+$.

Now (in $V$),  we claim that
$\exists \alpha < \lambda^+ \, \forall z \in S\, [ K_\alpha \notin z]$
(so that $\Gamma(x) \notin \cl(S)$).
If not, then we can fix  $W \in [\lambda^+]^{\lambda^+}$ and
$z \in S$ such that $K_\alpha \in z$ for all $\alpha \in W$.
Fix $\xi < \lambda$ such that $G_\xi \notin z$.
Since $\psi([W]^{< \omega}) = \lambda$, fix $s \in [W]^{< \omega}$
such that $\psi(s) = \xi$.
Say $s = \{\alpha_1,  \ldots , \alpha_n \}$.
Then $ G_\xi  = G_{\psi(s)} = H(s) \supseteq 
K_{\alpha_1} \cap \cdots \cap K_{\alpha_n} \in z$, a contradiction.
\end{proof}

\begin{proofof}{Lemma \ref{lemma-P-cover}}
Use Lemmas \ref{lemma-reflect} and \ref{lemma-reg},
with $\kappa = \lambda^+ = \aleph_n$.
\end{proofof}

In view of Lemma \ref{lemma-jon}, we can also prove Theorem \ref{thm-main}
replacing $\aleph_\omega$ with any other singular cardinal of cofinality
$\omega$,
since we can replace $\aleph_n$
in Lemma \ref{lemma-P-cover} by any successor to a regular cardinal.

\end{document}